\documentclass[letterpaper, 10 pt, conference]{ieeeconf}  % Comment this line out
\IEEEoverridecommandlockouts                              % This command is only
\usepackage[T1]{fontenc}
\usepackage[latin9]{inputenc}
\usepackage{amsmath}
\usepackage{amssymb}
\usepackage{graphicx}
\usepackage{dsfont}
\usepackage{mathdots}
\makeatletter

\pdfpageheight\paperheight
\pdfpagewidth\paperwidth

\@ifundefined{showcaptionsetup}{}{%
 \PassOptionsToPackage{caption=false}{subfig}}
\usepackage{subfig}
\makeatother


\newcommand{\spa}[1]{\text{span}\left(#1\right)}

\newcommand{\Proj}[2]{\mathcal{P}_{#1}\left(#2\right)}

\newcommand{\R}{\mathbb{R}}
\def\QED{~\rule[-1pt]{5pt}{5pt}\par\medskip}

\newtheorem{rem}{Remark}

\begin{document}

\title{Low-Rank and Low-Order Decompositions for Local System Identification}

\author{Nikolai Matni and Anders Rantzer
\thanks{N. Matni is with the Department of Control and Dynamical Systems, California Institute of Technology, Pasadena, CA.
 \tt{\small nmatni@caltech.edu}.}
 \thanks{A. Rantzer with Automatic Control LTH, Lund University, Box 118,
SE-221 00 Lund, Sweden. \tt{\small rantzer@control.lth.se}.}
\thanks{A. Rantzer gratefully acknowledges support of the LCCC Linnaeus Center and the eLLIIT Excellence Center at Lund University.}
\thanks{N. Matni was in part supported by NSF, AFOSR, ARPA-E, and the Institute for Collaborative Biotechnologies through grant W911NF-09-0001 from the U.S. Army Research Office. The content does not necessarily reflect the position or the policy of the Government, and no official endorsement should be inferred.}}% <-this % stops a space}
\maketitle
\begin{abstract}
As distributed systems increase in size, the need for scalable algorithms becomes more and more important.  We argue that in the context of system identification, an essential building block of any scalable algorithm is the ability to estimate local dynamics within a large interconnected system.  We show that in what we term the ``full interconnection measurement'' setting, this task is easily solved using existing system identification methods. We also propose a promising heuristic for the ``hidden interconnection measurement'' case, in which contributions to local measurements from both local and global dynamics need to be separated.  Inspired by the machine learning literature, and in particular by convex approaches to rank minimization and matrix decomposition, we exploit the fact that the transfer function of the local dynamics is low-order, but full-rank, while the transfer function of the global dynamics is high-order, but low-rank, to formulate this separation task as a nuclear norm minimization. 
\end{abstract}

\section{Introduction}
\label{sec:intro}
The new smart grid, the internet, and automated highway systems: all of these systems are characterized by their large scale, their distributed nature, and the sparse structure of their physical interconnections.  As these systems scale to larger and larger size, so too must the algorithms used to analyze and design them: thus local algorithms yielding global results become essential.  In general, such algorithms are not guaranteed to exist -- however, when additional structure is imposed on the system, it has been shown that there is indeed hope.

In the area of linear controller synthesis, for example, distributed systems systems with chordal structure \cite{RantzerChordal}, systems with favorable communication structures \cite{MickeyACC} and positive systems \cite{RantzerPositive} have all been shown to admit localized, and hence scalable, synthesis algorithms, with guaranteed global performance or stability guarantees.  Of course, none of these algorithms can be applied without first identifying the state-space parameters of the underlying large-scale distributed system.  Traditional system identification techniques such as subspace identification or prediction error are not computationally scalable -- furthermore, the former technique also destroys, rather than leverages, any \emph{a priori} information about the system's interconnection structure.  

We are not the first to make this observation, and indeed \cite{SIMsLargeScale} presents a local, structure preserving subspace identification algorithm for large scale (multi) banded systems (such as those that arise from the linearization of 2D and 3D partial differential equations), based on identifying local sub-system dynamics.  Their approach is to approximate neighboring sub-systems' states with linear combinations of inputs and outputs collected from a local neighborhood of sub-systems, and they show that the size of this neighborhood is dependent on the conditioning of the so-called structured observability matrix of the global system.

In this paper, we focus on the local identification problem, and leave the task of identifying the proper interconnection of these subsystems to future work, although we are also able to solve this problem in what we term the ``full interconnection measurement'' setting (to be formally defined in Section \ref{sec:problem}).  Our method is different from the approach suggested in \cite{SIMsLargeScale} in three respects: (1) we focus on identifying impulse response elements, rather than reconstructing state sequences, and (2) our methods are purely local, in that we do not require the exchange of information with any neighboring subsystems, and finally, (3) we do not need to assume a (multi) banded structure.  In light of this, we view our contribution as complementary to those presented in \cite{SIMsLargeScale}, and it will be interesting to to see if the two approaches can be combined in future work.

\begin{figure*}[!t]
\hfill{}\subfloat[Full interconnection measurements.]{\includegraphics[width=0.35\textwidth]{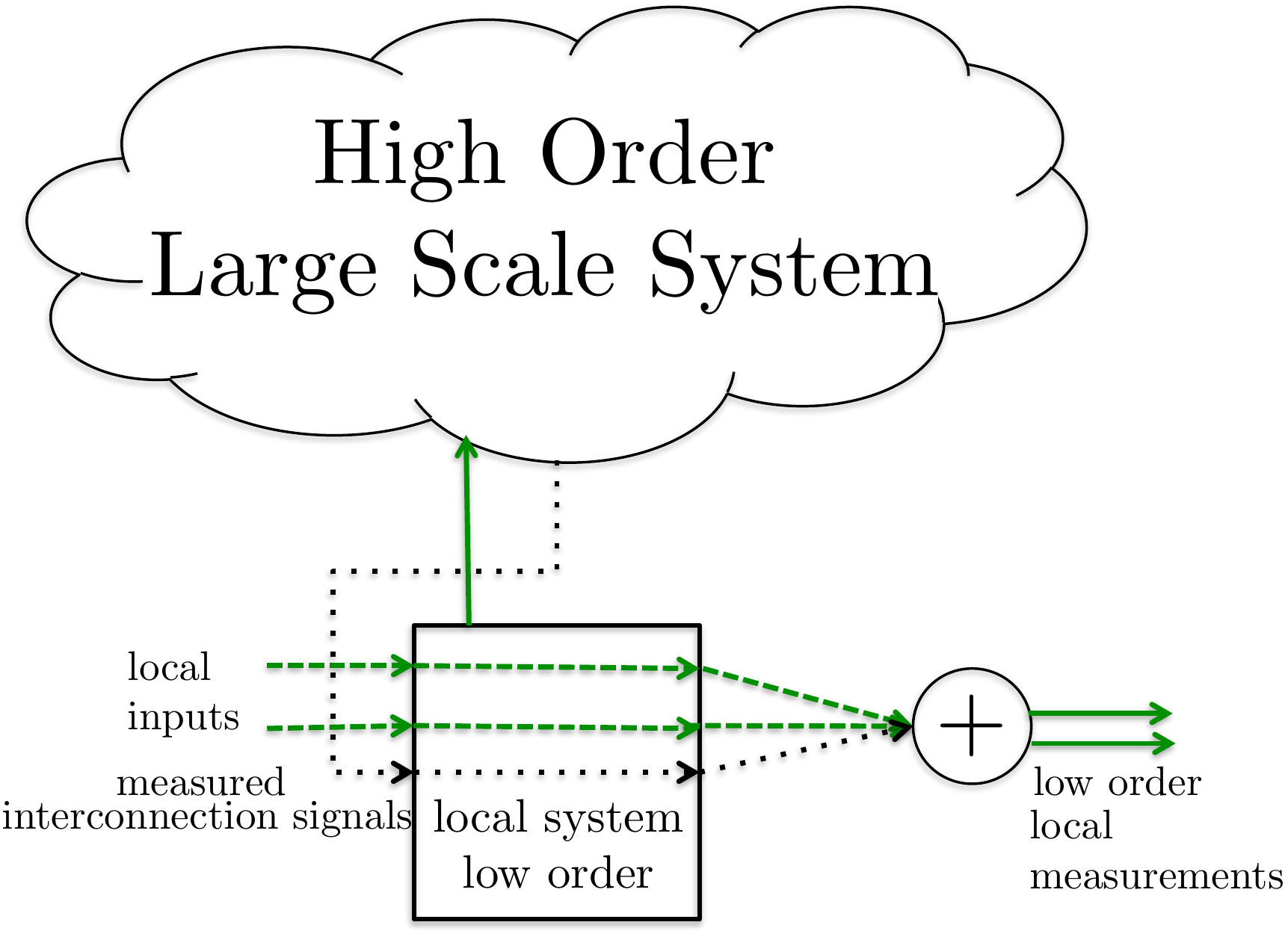}
}\hfill{}\subfloat[Hidden interconnection measurements]{\includegraphics[width=0.35\textwidth]{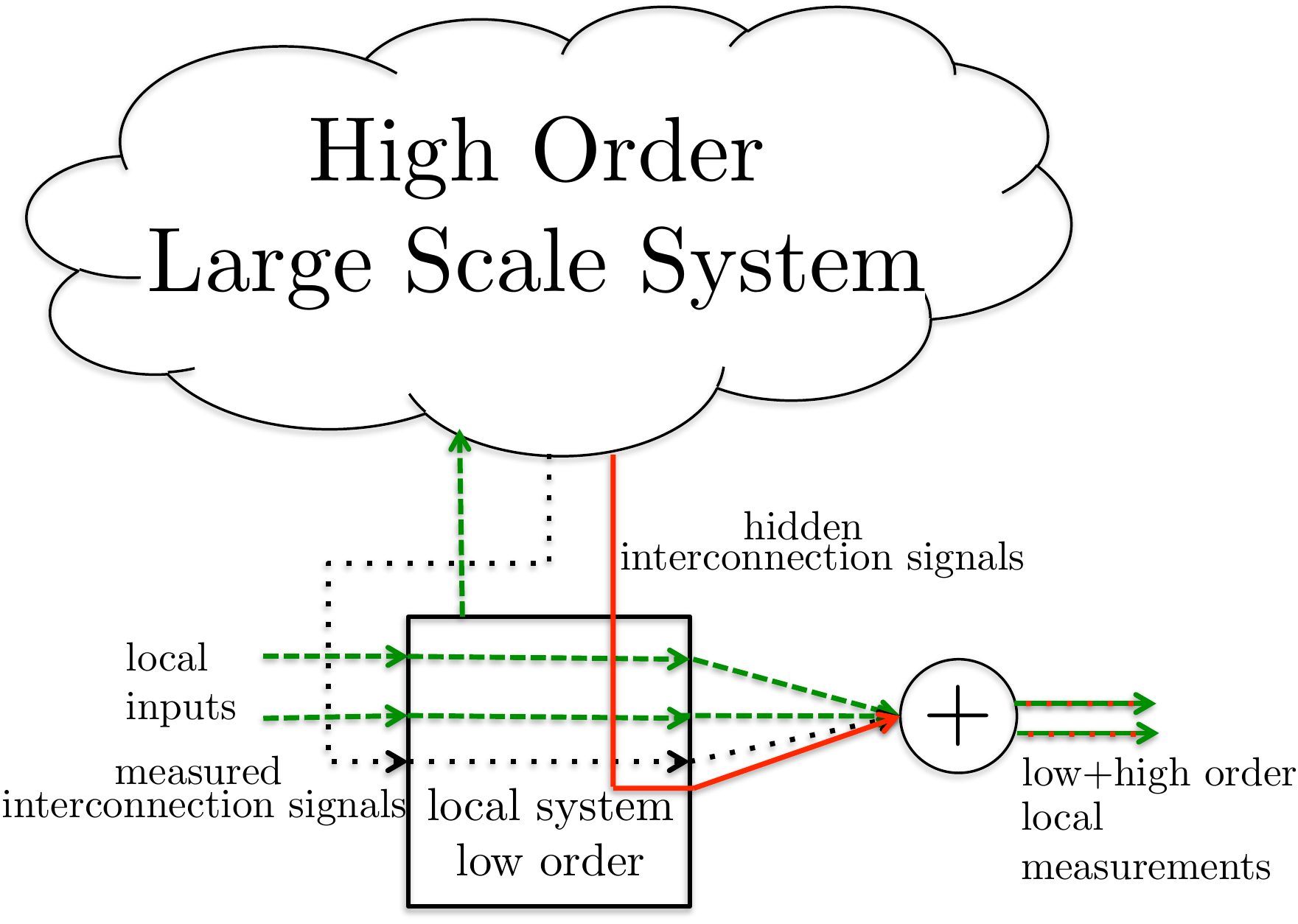}
}\hfill{}
\caption{\footnotesize Illustrated in Figures \ref{fig:local_vs_global} (a) and (b) are the full and hidden interconnection measurement cases, respectively.  Dashed green lines correspond to low-order signals, and dotted/solid black/red lines correspond to measured/hidden high-order interconnection signals.  In the full measurement case, the high order dynamics of the large scale system are isolated from the local measurements, as the interconnection signals can simply be treated as inputs to the system.  In the hidden interconnection measurement setting, high order global signals ``leak'' into our local measurements via the hidden interconnection signal (solid red), but do so through a low-rank transfer function.}
\label{fig:local_vs_global}
\end{figure*}

Our approach is based on two simple observations.  First, if all of the signals connecting the local sub-system to the global system, or \emph{interconnection signals}, can be measured, then under mild technical assumptions, the local observations are sufficient to identify both the local dynamics, and the coupling with the global system.  In effect, measuring the interconnection signals isolates the local sub-system, reducing the problem to a classical system identification problem.  Second, if an interconnection signal is not measured, then we have that the transfer function from local inputs and observed interconnection signals to local measurements naturally decomposes as the sum of two elements: one corresponding to local dynamics, which in general we expect to have \emph{full-rank, but low order}, and one corresponding to global dynamics, which will be of \emph{low-rank, but high order} (see Figure \ref{fig:local_vs_global} for a pictorial representation of both settings).  

Inspired by convex approaches to rank \cite{Fazel} and atomic norm minimization \cite{PariAtomicSysID} in system identification, and to matrix decomposition in latent variable identification in graphical models \cite{CPW10}, we conjecture that this difference in structure provides sufficient \emph{incoherence} (c.f. \cite{ChandraInco} and \cite{TroppInco} for examples of incoherence conditions) to allow the two signals to be separated through convex methods, in particular using nuclear norm minimization techniques.  Indeed a similar idea has been applied successfully to blind source separation problems \cite{TanakaBSS}.

This paper is organized as follows: in Section \ref{sec:problem}, we establish notation, and formally define the two variants of the problem to be solved, namely full and hidden interconnection measurement problems.  In Sections \ref{sec:full} and \ref{sec:hidden}, we provide nuclear norm minimization based algorithms for identifying local subsystem dynamics in both the full and hidden interconnection measurement settings, respectively.  We present numerical experiments supporting our approach in Section \ref{sec:numerical}, and end with conclusions and directions for future work in Section \ref{sec:conclusion}.

\section{Problem Formulation}
\label{sec:problem}
\subsection{Notation}
For a matrix 
\begin{equation}
X = \begin{bmatrix} X_0 & X_1 \dots & X_{2N} \end{bmatrix}
\end{equation}
we define the Hankel operator $\mathcal{H}(X)$ to be
\begin{equation}
\mathcal{H}(X) := \begin{bmatrix} X_1 & X_2 & \dots &  X_{N} \\
					           X_2 & X_3 & \iddots &  X_X \\
					           \vdots & \iddots & \ddots &  \vdots \\
					           X_{N} & X_{N+1} & \dots & X_{2N} \end{bmatrix},
\end{equation}
and its Fourier transform to be given by
\begin{equation}
\mathcal{F}(X)(e^{j\omega_k}) = 
\sum_{t=0}^{2N-1} X_t e^{-{j\omega_kt}}
\end{equation}
for $\omega_k = \frac{\pi k}{N}$, $k\in \{0,\dots,2N-1\}$.

For a set of measurements $\{m^i_t\}_{t=0}^N$, $m^i_t \in R^C$, and natural numbers, $N$, $M$ and $r$, with $N$ even, we define $M^i_{N,M,r} \in \R^{C(r+1)\times(M+1)}$ by
\begin{equation}
M^i_{N,M,r} := \begin{bmatrix}
m^i_{N-M} & m^i_{N-(M-1)} & \dots & m^i_N \\
m^i_{N-(M+1)} & m^i_{N-M}  & \dots & m^i_{N-1} \\
\vdots & \vdots & \ddots & \vdots \\
m^i_{N-(M+r)} & m^i_{N-(M-1+r)} & \dots & m^i_{N-r} 
\end{bmatrix},
\label{eq:NMR}
\end{equation}
where we adopt the convention that $m^i_t = 0$ for all $t<0$.  When $N$, $M$ and $r$ are clear from context, we drop the subscripts and simply denote the matrix by $M^i$.  

For a general matrix $M$, we let $\|M\|_F$ denote its Froebenius norm, i.e. $\|M\|_F^2 = \mathrm{trace}M^\top M$, and $\|M\|_*$ denote its nuclear norm, i.e. $\|M\|_* = \sum_i \sigma_i$, where $\sigma_i$ are the singular values of $M$.

For a subspace $\mathcal{S}$, we denote by $\Proj{\mathcal{S}}{\cdot}$ the orthogonal projection operator onto $\mathcal{S}$ with respect to the euclidean inner-product, and by $\mathcal{S}^\perp$ the orthogonal complement of the subspace, once again with respect to the euclidean inner-product.

\subsection{Distributed systems with sparse interconnections}

We consider a distributed system comprised of $n$ linear time
invariant (LTI) sub-systems, which interact with each other according
to a physical interaction graph
$\mathcal{G}=(\mathcal{X},\, E)$. We denote by $i\in\mathcal{X}$ the $i^{\textnormal{th}}$
node in the graph, and by $x^{i}$ the state of the corresponding
sub-system. We assume that each subsystem $i\in\mathcal{X}$ has its own
control input $u^{i}$ and centered white noise process noise $w^{i}$ (satisfying  $\mathbb{E}[w_t^{i}{w_t^{j}}^{\top}]= W^{ij}$, $\mathbb{E}[w^{i}_s{w^{j}_t}^{\top}]=0\ \forall s \neq t$),
and that plants physically interact with each other according to $E$.
In particular, an edge $e_{ij} \in E$ is non-zero
if and only if subsystem $j$ directly affects the dynamics of subsystem $i$.  Defining the neighbor set of node $i$ as $\mathcal{N}_i = \{ j \in \mathcal{X} \, : \, e_{ij}\neq 0\}$, we can then write the dynamics of each subsystem as
\begin{equation}
x^i_{t+1} = A^{ii}x^i_t + \sum_{j \in \mathcal{N}_i} A^{ij}x^j_t + B^iu^i_t + w^i_t,
\label{eq:subsystem}
\end{equation}
with initial conditions $x_i(0)=0$, subsystem state $x^i \in \R^{n_i}$, neighboring subsystem states $x^j_t \in \R^{n_j}$, subsystem input $u^i \in \R^{p_i}$ and subsystem process noise $w^i_t \in \R^{n_i}$.  For reasons that will become apparent, we will refer to the signals $(A^{ij}x^j_t)_{t=0}^N$ as the interconnection signals at node $i$ over a horizon $N\geq 0$.

\subsection{Local and interconnection observations}

In the following we distinguish between two types of observations that can be collected at node $i$.  The first, which we call local observations, correspond to standard measurements of the local state, i.e. we call $y^i_t \in \R^{q_i}$, as given by
\begin{equation}
y^i_t = C^i x^i_t + D^i u^i_t +  \nu^i_t,
\label{eq:local_obs}
\end{equation}
the local state observations at time $t$, with $\nu^i_t \in \R^{q^i}$ a centered white noise process.

The second, which we term interconnection observations, correspond to measurements of incoming signals from neighboring nodes, i.e. we call $z^i_t \in \R^{m_i}$, as given by
\begin{equation}
\begin{array}{rcl}
z^i_t &=& \bar{C}^i \bar{x}^i_t + \bar{\nu}^i_t,
\end{array}
\label{eq:observations}
\end{equation}
the interconnection observations at time $t$, where $\bar{x}^i_t = (x^j_t)_{j\in\mathcal{N}_i} \in \R^{\sum_{j\in\mathcal{N}_i} n_j},$ and $\bar{\nu}_t^i \in \R^{\sum_{j\in\mathcal{N}_i} n_j}$ is a centered white noise process.

\subsection{Local system identification}

Our system identification goal is to identify, up to a similarity transformation, the the tuple $(A^{ii}, B^{i}, C^{i}, D^{i})$ given only the time history of $(u^i, y^i, z^i)$ -- that is to say we seek a local estimation procedure for the subsystem dynamics.  This task is non-trivial as the subsystem is connected to the remaining full system, and thus even identifying the true order of the local subsystem can be challenging. 

In the sequel, we assume that the full system is Hurwitz, that $(A^{ii}, C^i)$ is observable, and without loss that each $C^i$ has full row rank, and once again distinguish between two cases.  The first is when we have that all interconnection signals are contained within the linear span of the interconnection observations -- we refer to this case as the \emph{full interconnection measurement} case.  Formally, this can be stated as 
\begin{equation}
A^{ij}x^j_t \in \spa{\bar{C}^i \bar{x}^i_t}, \, \forall t\geq 0, \, \forall j \in \mathcal{N}_i, 
\label{eq:cutoff}
\end{equation}
or more succinctly, that there exists a linear transformation $\mathcal{L}_{ij}$ such that
\begin{equation}
A^{ij} = \mathcal{L}_{ij}(\bar{C}^i), \, \forall j \in \mathcal{N}_i.
\label{eq:cutoff2}
\end{equation}
We also define $\mathcal{L}_i$ as the linear operator 
\begin{equation}
\mathcal{L}_i := [\mathcal{L}_{ij_1},\dots,\mathcal{L}_{ij_{|\mathcal{N}_i}|}]
\end{equation}
such that
\begin{equation}
\sum_{j\in\mathcal{N}_i}A^{ij}x^j_t = \mathcal{L}_i(\bar{C})\bar{x}^i_t, \ \forall t\geq 0.
\end{equation}

We will show that under mild coordination with neighboring subsystems, we are able to identify $(A^{ii}, B^{i}, C^{i}, D^{i})$  (to within the accuracy allowable by the noise) using only local information.  Intuitively, by measuring these connecting signals, they can be treated as inputs to the subsystem, effectively isolating node $i$ from the global dynamics (see Figure \ref{fig:local_vs_global}(a)) -- however, in order to ensure persistence of excitation under this setting, non-local elements of randomness need to be injected into the system, hence the need for coordination. 

The second case, which we call the \emph{hidden interconnection measurement} setting, occurs when not all interconnection signals are observed, i.e. when conditions \eqref{eq:cutoff} or \eqref{eq:cutoff2} do not hold.  The local dynamics can no longer be isolated from the global dynamics due to these unobserved interconnection signals -- as such, our full interconnection measurement method would lead to the identification of a high order local model due to the ``hidden'' connection to the full system (see Figure \ref{fig:local_vs_global}(b)).  Inspired by the success of convex methods for sparse and low-rank decomposition techniques in identifying latent variables in graphical models \cite{CPW10}, and for blind source separation \cite{TanakaBSS}, we propose a convex programming method for identifying and separating out the local low-order dynamics from the global high-order dynamics, which are due to the hidden connection with the full system.

\section{Full inter-connection measurements}
\label{sec:full}
We begin by assuming that \eqref{eq:cutoff} and \eqref{eq:cutoff2} hold, and consider the case when all noise terms are identically zero.  A robust variant of our solution will be presented at the end of this section when noise is present in the system.

For any $t\geq 0$, we may then write
\begin{equation}
y^i_t = \sum_{k=0}^t s^i_k \begin{bmatrix} u^i_{t-k} \\ z^i_{t-k} \end{bmatrix},
\label{eq:impulse}
\end{equation}
with $s^i_0 = [D^i,0]$, $s^i_t = C^i (A^{ii})^{t-1} [B^i,\mathcal{L}(\bar{C}^i)]$ the subsystem's impulse response elements.

With this in mind, fix natural numbers $N$, $M$ and $r$, with $N$ even, and let $v^i_t = [{u^i_t}^\top, \, {z^i_t}^\top]^\top$, $V^i_{N,M,r}$ be given by \eqref{eq:NMR}, and
%\begin{equation}
%V^i = \begin{bmatrix}
%v^i_{N-M} & v^i_{N-(M-1)} & \dots & v^i_N \\
%v^i_{N-(M+1)} & v^i_{N-M}  & \dots & v^i_{N-1} \\
%\vdots & \vdots & \ddots & \vdots \\
%v^i_{N-(M+r)} & v^i_{N-(M-1+r)} & \dots & v^i_{N-r} 
%\end{bmatrix} %\in \R^{(r+1)p_i \times (M+1)}
%\label{eq:Vi}
%\end{equation}
%\begin{equation}
%U^i = \begin{bmatrix}
%u^i_{N-M} & u^i_{N-(M-1)} & \dots & u^i_N \\
%u^i_{N-(M+1)} & u^i_{N-M}  & \dots & u^i_{N-1} \\
%\vdots & \vdots & \ddots & \vdots \\
%u^i_{N-(M+r)} & u^i_{N-(M-1+r)} & \dots & u^i_{N-r} 
%\end{bmatrix} %\in \R^{(r+1)p_i \times (M+1)}
%\label{eq:Ui}
%\end{equation}
%\begin{equation}
%Z^i = \begin{bmatrix}
%z^i_{N-M} & z^i_{N-(M-1)} & \dots & z^i_N \\
%z^i_{N-(M+1)} & z^i_{N-M}  & \dots & z^i_{N-1} \\
%\vdots & \vdots & \ddots & \vdots \\
%z^i_{N-(M+r)} & z^i_{N-(M-1+r)} & \dots & z^i_{N-r} 
%\end{bmatrix} %\in \R^{(r+1)m_i \times (M+1)}
%\label{eq:Zi}
%\end{equation}
\begin{equation}
Y^i = \begin{bmatrix} y^i_{N-M} & y^i_{N-(M-1)} & \dots & y^i_N \end{bmatrix}
\label{eq:Yi}
\end{equation}
\begin{equation}
S^i = \begin{bmatrix} s^i_0 & s^i_1 & \dots & s^i_r \end{bmatrix}.
\label{eq:Si}
\end{equation}

Choosing $r = N$, we may then write
\begin{equation}
Y^i = S^i V^i.
\label{eq:consistency}
\end{equation}
 Thus we seek conditions under which \eqref{eq:consistency} has a unique solution -- i.e. we seek conditions under which \[ V^i \in \R^{(N+1)(p_i+m_i)\times(M+1)}\] has a right inverse, yielding the solution
 \begin{equation}
S^i = Y^i (V^i)^\dag,
 \end{equation}
 where $X^\dag$ denotes the pseudo-inverse of $X$.
 
   A necessary condition, that we assume holds in the sequel, is that $M$ is sufficiently large such that $M+1 \geq (N+1)(p_i+m_i)$.
 
 \begin{rem}
 One may choose to approximate outputs as coming from a finite impulse response system of order $r$ by choosing $r<N$; as the system is assumed to be stable, picking a sufficiently large $r$ then allows for a computational gain without sacrificing accuracy.  In this case, the aforementioned necessary condition then becomes $M+1 \geq (r+1)(p_i+m_i)$.
 \end{rem}
 
Next we characterize necessary and sufficient conditions for $V^i$ to have full row-rank. In order to make the analysis more transparent, introduce the auxiliary matrices $U^i$ and $Z^i$, constructed from $\{u^i_t\}_{t=0}^N$ and $\{z^i_t\}_{t=0}^N$, respectively, 
%\begin{equation}
%U^i = \begin{bmatrix}
%u^i_{N-M} & u^i_{N-(M-1)} & \dots & u^i_N \\
%u^i_{N-(M+1)} & u^i_{N-M}  & \dots & u^i_{N-1} \\
%\vdots & \vdots & \ddots & \vdots \\
%u^i_{N-(M+r)} & u^i_{N-(M-1+r)} & \dots & u^i_{N-r} 
%\end{bmatrix} %\in \R^{(r+1)p_i \times (M+1)}
%\label{eq:Ui}
%\end{equation}
%\begin{equation}
%Z^i = \begin{bmatrix}
%z^i_{N-M} & z^i_{N-(M-1)} & \dots & z^i_N \\
%z^i_{N-(M+1)} & z^i_{N-M}  & \dots & z^i_{N-1} \\
%\vdots & \vdots & \ddots & \vdots \\
%z^i_{N-(M+r)} & z^i_{N-(M-1+r)} & \dots & z^i_{N-r} 
%\end{bmatrix} %\in \R^{(r+1)m_i \times (M+1)}
%\label{eq:Zi}
%\end{equation}
 and note that \[\text{rank}(V^i) = \text{rank}\left(\begin{bmatrix} U^i \\ Z^i \end{bmatrix}\right).\]
 
Therefore, necessary and sufficient conditions are that each of (i) $U^i$ and (ii) $\Proj{{U^i}^\perp}{Z^i}$ (the projection of $Z^i$ onto the orthogonal complement of the row space of $U^i$) have full row rank.  Condition (i) is easily satisfied (with probability one) by choosing $u^i_t$ to be a white random process -- we therefore assume this holds and focus on condition (ii).
 
 It should be immediate to see that if no other inputs are administered to the system then $\Proj{{U^i}^\perp}{Z^i}=0$, as the system's trajectory lies entirely in the span of the row space of $U^i$.  Therefore, let $\mathcal{A}_i := \{j \in \mathcal{X} \, : \, u^j \not \equiv 0 \}$ denote the set of ``active'' inputs in the rest of the system, and let $u^{-i}_t=\left( u^j_t \right)_{j\neq i \in \mathcal{A}_i}$.  
 
 Then $Z^i \in \spa{U^i, U^{-i} }$, where $U^{-i} $ is generated by $\{ u^{-i}_t \} _{t=0}^N$.  If (i) the transfer function from $u^{-i}$ to $z^i$ has full row rank, and (ii) sufficiently many active inputs are present (specifically, a number greater than or equal to $m_i$), and chosen to be such that $U^{-i}$ is full row rank (which, again, is generically true for white input processes), then indeed $\Proj{{U^i}^\perp}{Z^i}$ will have full row rank.
 
Thus we see that through a marginal amount of coordination (signaling other subsystems to inject exciting inputs into the system), a purely local estimation procedure can be used to exactly recover the first $N$ impulse response elements $s_0, \, \dots, \, s_N$ of the local subsystem, to which standard realization procedures can then be applied to extract (up to a similarity transformation), the tuple $(A^{ii}, [B^i,\mathcal{L}(\bar{C}^i)], C^i, D^i)$. 

\subsection{A robust variant}

Following \cite{Fazel}, we can formulate a robust variant of our previous approach when the noise terms are non-zero.  Defining 
\begin{equation}
\Delta^i := Y^i  - S^iV^i
\end{equation}
we then solve the following nuclear norm minimization
\begin{equation}
\begin{array}{rl}
\text{minimize}_{S^i} & \|\mathcal{H}(S^i)\|_* \\
\text{s.t.} & \|\Delta^i\|_F \leq \delta
\end{array}
\label{eq:robust}
\end{equation}
where $\delta$ is a tuning parameter that ensures consistency of the estimated impulse response elements with the observed data.  Note that this approach can also be suitably modified to accommodate bounded noise \cite{Fazel}, or unbounded noise with known covariance \cite{LjungNewOld}, or to handle missing time points in the output signal data as described in \cite{LHV13}.

\section{Hidden inter-connection measurements}
\label{sec:hidden}
When condition \eqref{eq:cutoff} does not hold, the local identification task becomes much more difficult -- by not measuring all of the connecting signals, global high-order dynamics ``leak'' into our local estimation procedure (see Figure \ref{fig:local_vs_global}(b)).  Inspired by sparse and low-rank decomposition methods used to identify latent variables in graphical models \cite{CPW10}, and by Hankel rank minimization techniques used in blind source separation problems \cite{TanakaBSS}, this section proposes a regularized variant of program \eqref{eq:robust} that has shown promise in numerical experiments. 

Formal results proving the success of this technique (analogous to those found in \cite{CPW10,RFP10}) are the subject of current work.  This subsection aims rather to provide some intuition and justification for the method.  In particular, define the number of hidden signals at node $i$ to be
\begin{equation}
k_i = \sum_{j\in\mathcal{N}_i} \dim\left( \Proj{\spa{\bar{C}}^\perp}{\spa{A^{ij}}}\right)
%k_i = \min_{\mathcal{L}_i}\left\{\sum_{j\in\mathcal{N}_i} \dim \Proj{\spa{\mathcal{L}_i(\bar{C}^i)}^\perp}{\spa{A^{ij}}}\right\},
\end{equation}
that is to say, the dimension of the subspace of the hidden interconnection signals.  

We may then write, analogous to \eqref{eq:impulse}

\begin{equation}
y^i_t = \sum_{k=0}^t s^i_k \begin{bmatrix} u^i_{t-k} \\ z^i_{t-k} \end{bmatrix} 
+ \sum_{k=0}^t h^i_k \begin{bmatrix} u^i_{t-k} \\ z^i_{t-k} \\ u^{-i}_{t-k} \end{bmatrix} .
\end{equation}
where the $s^i$ are once again the impulse response elements of the local-subsystem, whereas $h^i_0 = 0$,  
%\begin{multline}
%h^i_t =\\
%   \indent \begin{bmatrix} 0 & \dots & C^i & \dots & 0 \end{bmatrix} A^{t-1} \begin{bmatrix} 0 & 0 & 0 \\ \vdots & \vdots & \vdots \\ B^i & \mathcal{L}_i(\bar{C}^i) & B^{\mathcal{A}_i} \\ \vdots & \vdots & \vdots \\ 0 & 0  & 0 \end{bmatrix}  \\  \indent \indent - \begin{bmatrix} s^i_t & 0 \end{bmatrix}
%\end{multline}
and $(h^i_t)$ are the impulse response elements describing the global dynamics that are ``leaking'' in to our subsystem via the hidden interconnection signals.  Let $w^i_t = [ {v^i_t}^\top, \,{u^{-i}_t}^\top]^\top$, and  $W^i$ be as in \eqref{eq:NMR}, and
\begin{equation}
H^i = \begin{bmatrix} H^i_{N-M} & H^i_{N-(M-1)} & \dots & H^i_N \end{bmatrix},
\end{equation}
allowing us to write
\begin{equation}
Y^i = \begin{bmatrix} S^i & H^i \end{bmatrix} \begin{bmatrix} V^i \\ W^i \end{bmatrix}.
\label{eq:consistency2}
\end{equation}

We now make the key observation that the transfer function $H(e^{j\omega_k}) = \mathcal{F}(H^i)$ can have rank at most $k_i$, the number of hidden interconnection signals.  In all of the following, we assume that the transfer function from $(u^i,z^i)$ to $y^i$ is full rank, and that
\begin{equation}
\min\left(p_i + m_i, \, q_i\right) > k_i
\end{equation}
holds.  Specifically, we ask that both the dimension  $q_i$ of the subspace spanned by our local observations, and the dimension $p_i + m_i$ of the subspace spanned by the ``inputs'' $u^i$ and $z^i$, be larger than the dimension $k_i$ of the subspace spanned by the hidden interconnection signals.  Interpreted in terms of the rank of transfer functions, we ask that the rank of the local component of transfer function from $(u^i,z^i)$ to $y^i$, given by $\min\left(p_i + m_i, \, q_i\right)$ under our full rank assumption, be larger than the rank $k_i$ of the global component of the transfer function from $(u^i,z^i)$ to $y_i$.

If these conditions hold, we then have a \emph{structural} means of distinguishing between the two components of the impulse response of the local subsystem.  First, we expect $\mathcal{H}(S^i)$ to have low rank, as it describes the low-order dynamics of the local model, whereas $\mathcal{H}(H^i)$ will not, as it corresponds to the high-order global dynamics that leak in via the hidden connecting signals.  Secondly, by our local full rank assumption, we will have that $S\left(e^{j\omega_k}\right)$ is full rank (with rank $\min\left(p_i + m_i, \, q_i\right)$), whereas  $\text{rank}\left( H\left(e^{j\omega_k}\right)\right) \leq k_i$; as mentioned above, the hidden interconnection signals act as a structural ``choke'' point, limiting the rank of the interconnecting transfer function.  

This suggests a natural decomposition of the impulse response elements of $y^i$ into a local full rank but low-order component with simple dynamics, and a hidden high-order but low-rank component.  Using the nuclear-norm heuristic for low-rank approximations \cite{RFP10}, we may then modify program \eqref{eq:robust} to control the rank of $H\left(e^{j\omega_k}\right)$:
\begin{equation}
\begin{array}{rl}
\text{minimize}_{S^i, \, H^i} & \|\mathcal{H}(S^i)\|_* \\
\text{s.t.} & \Delta^i = 0  \\
& \|H\left(e^{j\omega_k}\right)\|_* \leq \delta_h, \\ & \indent \indent \indent \omega_k = \frac{2\pi k}{M}, \, k=0,\,1,\dots,M-1
\end{array}
\label{eq:hidden}
\end{equation}
where now
\begin{equation}
\Delta^i = Y^i - \left( \begin{bmatrix} S^i & H^i \end{bmatrix} \begin{bmatrix} V^i \\ W^i \end{bmatrix}\right),
\end{equation}
and $\delta_h$ is an additional tuning parameter used to control the rank of $H\left(e^{j\omega_k}\right)$ across frequencies.   When noise is present, we relax the constraint on $\Delta^i$ to $\|\Delta^i\|_F \leq \delta$, as in the robust variant of the full interconnection measurement case.

This method is, however, non-local in that it requires the communication of $U^{-i}$ to node $i$ in order to implement it.  In light of this, we also suggest the following local approximation to \eqref{eq:hidden}.  In particular, we define 
\begin{equation} \tilde{\Delta}_i = Y^i - \left( \begin{bmatrix} S^i & H^i \end{bmatrix} \begin{bmatrix} V^i \\ V^i \end{bmatrix}\right)
\end{equation}
and propose solving
\begin{equation}
\begin{array}{rl}
\text{minimize}_{S^i, \, H^i} & \|\mathcal{H}(S^i)\|_* \\
\text{s.t.} & \| \tilde{\Delta}^i \|_F \leq \delta  \\
& \|H\left(e^{j\omega_k}\right)\|_* \leq \delta_h, \\ & \indent \indent \indent \omega_k = \frac{2\pi k}{M}, \, k=0,\,1,\dots,M-1.
\end{array}
\label{eq:hidden2}
\end{equation}

Essentially, we treat the unknown active inputs $U^{-i}$  as disturbances entering the system through $H^i(e^{j\omega_k})$, and therefore allow $\tilde{\Delta}_i$ to deviate from 0, but still insist on consistency with the observed data.

\section{Numerical Experiments}
\label{sec:numerical}
\begin{figure}
\centerline{\includegraphics{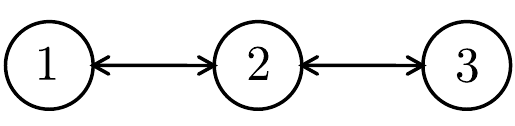}}
\caption{\footnotesize The graph depicts the physical interconnection structure of the three-subsystem chain.}
\label{fig:chain}
\end{figure}
We consider the following three subsystem chain (as illustrated in Figure \ref{fig:chain}), with $x_t, \, w_t \in \R^9$ and $u \in \R^5$,
\begin{equation}
x_t = Ax_t + Bu_t + w_t
\end{equation}
with $A$ and $B$ given as in equations \eqref{eq:A} and \eqref{eq:B} (found at the end of the paper), and identically and independently distributed $w_t\sim \mathcal{N}(0,.01^2 I)$.  Each node has a state $x^i_t \in \R^3$, which we assume are ordered such that \[ x_t = \begin{bmatrix} x^1_t \\ x^2_t \\ x^3_t \end{bmatrix} .\]

We will consider the task of identifying node $1$'s system parameters, namely we seek to identify the tuple $(A^{11},B^1,C^1,D^1)$ where
\begin{equation}
 A^{11} = 
 \begin{bmatrix}
  0.2839    &0.2125 &  -0.3097\\
    0.1528&   -0.3525 &   0.2400\\
    0.0183 &  -0.1709 &  -0.0109
    \end{bmatrix}, \end{equation}
    
 \begin{equation}
  B^1 =   \begin{bmatrix}
    0.6394 &  -0.3201 \\
    0.8742  & -0.1374 \\
    1.7524   & 0.6158
\end{bmatrix}
\end{equation}
\begin{equation}
C^1 = \begin{bmatrix}
  0.6348   &-0.1760 &  -0.1274 \\
    0.8204   & 0.5625&    0.5542
    \end{bmatrix}
    \end{equation}
    \begin{equation}
    D^1 = 
    \begin{bmatrix}
    -1.0973  &  1.4047 \\
   -0.7313  & -0.6202
    \end{bmatrix}
\end{equation}
given local observations $y^1_t = C^1 x^1_t + \nu^1_t$, with $\nu^1_t\sim \mathcal{N}(0,.01^2 I)$ and varying amounts of interconnection measurements.  Note that in this system $\bar{x}^i_t = x^2_t$, and that indeed this fact remains true regardless of the number of subsystems in the chain. 

We begin with the full interconnection measurement setting, with measurement noise $\bar{\nu}^1_t\sim \mathcal{N}(0,.01^2 I)$ and
\begin{equation}
z^1_t = \begin{bmatrix}
0.4895  &  0.6449   & 0.4762 \\
   -1.5874  &  0.1367 &   0.6874 \\
    0.8908 &   0.1401 &   0.9721
\end{bmatrix}x^2_t + \bar{\nu}^1_t =: \bar{C}^i x^2_t + \bar{\nu}^1_t.
\end{equation}
It is easily verified that $\bar{C}^i$ is invertible, and thus satisfies \eqref{eq:cutoff} and \eqref{eq:cutoff2}.  Solving program \eqref{eq:robust} with $N=600$, $M = 300$, $r=21$ and $\delta = 0.5$, we obtain an estimation error of $\| \hat{S}^i - S^i \|_F = .008$, relative to $\|S\|_F = 2.871$; i.e. we recover the impulse response elements to within the limits set by the noise.  Additionally, $\text{rank}\left(\mathcal{H}(S^i)\right) = 3$, the true order of the system.

\begin{figure*}
\centerline{\includegraphics[width=\textwidth]{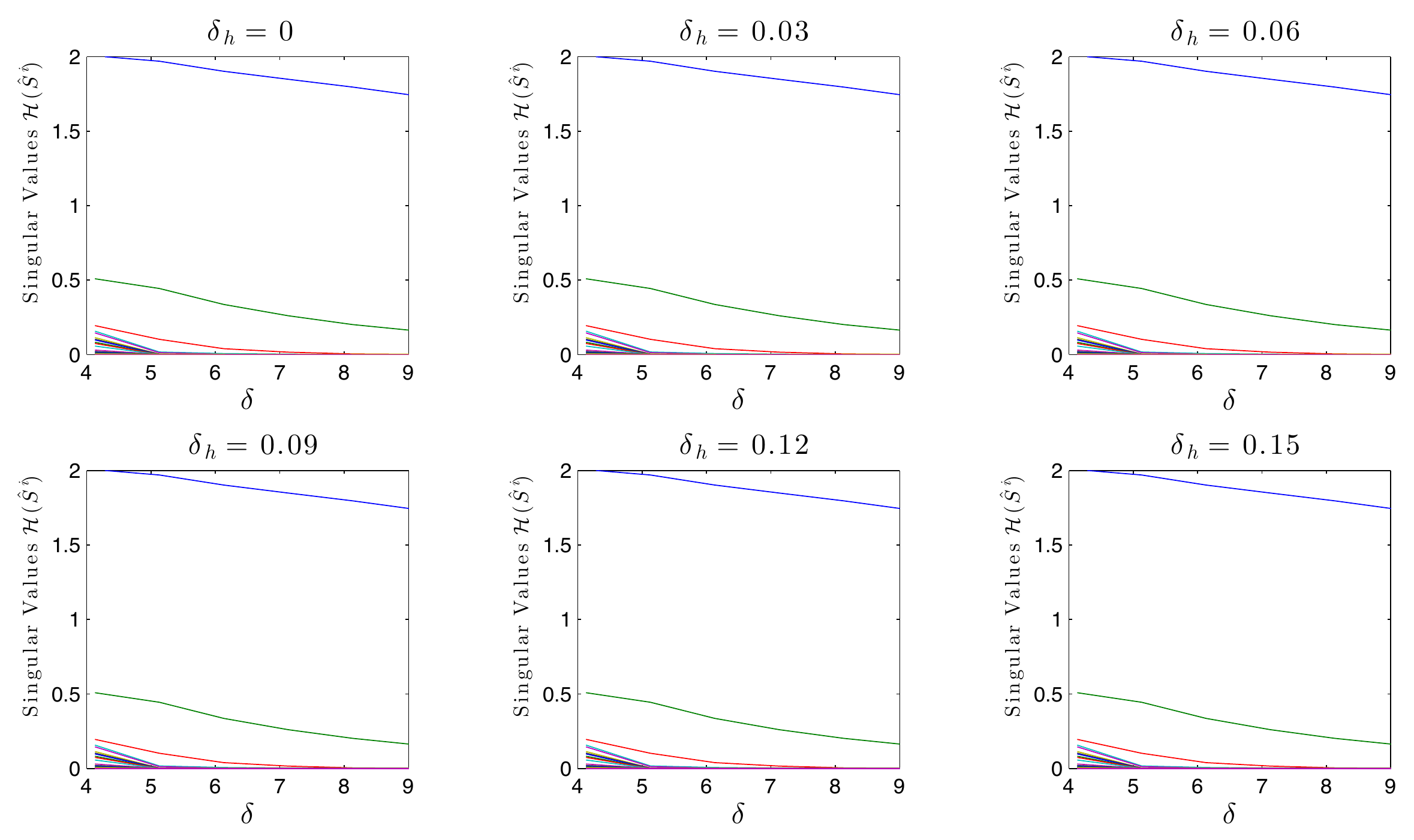}}
\caption{\footnotesize By examining how the values of the singular values of $\mathcal{H}(\hat{S}^i)$ vary across different values of $\delta$ and $\delta_h$, the order of the local sub-system is correctly identified as three.}
\label{fig:singularvals}
\end{figure*}

Next we consider the case where we have hidden interconnection signals.  In particular, we let  
\begin{equation}
z^1_t = \begin{bmatrix}
0.4895  &  0.6449   & 0.4762 \\
   -1.5874  &  0.1367 &   0.6874 \\
\end{bmatrix}x^2_t =: \bar{C}^i x^2_t.
\end{equation}
Once again, we easily verify that $\bar{C}^i$ has full row-rank of 2, and therefore conclude that the dimension of the hidden interconnection subspace is 1, which is less than $p_1+m_1 = 4$ and $q_1 = 2$. We solve program \eqref{eq:hidden2} with $N=600$, $M = 300$, $r=21$, $\delta_h = .05$ and $\delta = 4.5$, and obtain an estimation error of $\| \hat{S}^i - S^i \|_F = .093$, relative to $\|S\|_F = 2.871$; although our error is above the noise level, it is still a reasonable estimate of the local dynamics.  Most importantly we believe, however, is that (i) the top three singular values of $\mathcal{H}(S^i)$ were at least an order of magnitude larger than the remaining singular values for a fairly broad range of  $\delta$ and $\delta_h$ (see Figure \ref{fig:singularvals}), and that (ii) the rank of each $H(e^{j\omega_k})$ term was correctly identified as 1 for all values of $\delta_h \in [0,0.15]$ across a broad range of values of $\delta$.  Indeed, numerical experiments seem to suggest that the method is well suited to identifying the true order of the local dynamics, and the dimension of the hidden interconnection subspace, opening up the possibility of further refining results using parametric methods.

\section{Conclusion}
\label{sec:conclusion}
We presented a nuclear norm minimization based approach to separating local and global dynamics from local observations, and argued that this method can be used as part of a distributed system identification algorithm.  In particular, we noted that when all interconnection signals can be measured, the problem essentially reduces to a classical system identification problem.  When some interconnection signals are not measured, we exploit the fact that the transfer function from $(u^i,z^i)$ to $y^i$ naturally decomposes into a local contribution that is low-order, but full rank, and a global contribution that is high-order, but low rank to formulate the local system identification problem as a matrix decomposition problem amenable to convex programming. 

In future work, we will look to develop non-asymptotic consistency results for our estimation procedure, analogous to those found in \cite{CPW10,RFP10,PariAtomicSysID}.  It is also of importance to develop a principled method for interconnecting our local sub-systems properly to ultimately yield an accurate global model, analogous to the algorithm presented in \cite{SIMsLargeScale}.  Finally, more numerical experiments need to be conducted to further validate the efficacy of this method, especially on real world, as opposed to synthetic, data.

\section{Acknowledgements}
The authors would like to thank Lennart Ljung for the fruitful suggestion of formulating the problem as a Frobenius norm constrained optimization.

\bibliographystyle{/Users/nmatni/Documents/Publications/ACC13_duality/IEEEtran}
\bibliography{/Users/nmatni/Documents/Publications/biblio/comms,/Users/nmatni/Documents/Publications/biblio/decentralized,/Users/nmatni/Documents/Publications/biblio/matni,/Users/nmatni/Documents/Publications/biblio/regularization,/Users/nmatni/Documents/Publications/biblio/sys_id}

\begin{figure*}
\begin{equation}
\label{eq:A}
A = \begin{bmatrix}  
    0.2839  & 0.2125   & 	-0.3097   & 	0.1843 &   0.0775   &	-0.1358      &   0  &       0   &      0 \\
    0.1528 &  -0.3525   & 	0.2400    &	0.0976  & -0.1246   &	-0.0821       &  0        & 0        & 0\\
    0.0183  & -0.1709   & 	-0.0109   &	-0.3269  & -0.0005   & 0.1012       	 & 0 &        0  &       0 \\
    0.0857   & 0.3037   & 	-0.1947   &	0.0914    &0.3916   	& 0.3797   	& 0.0774 &  -0.0510  &  0.2253\\
   -0.1698  & -0.1557   & 	-0.1865   &	 0.2742   & 0.2066 	 & -0.5958 &	   0.3695&    0.1370   &-0.4422\\
    0.4134   & 0.1407    & 	0.2100    &	0.1776    &0.0653  	& -0.2677   &	 0.1827   &-0.2593   & 0.0085\\
         0        & 0         &	 0  		&	 -0.5795   &-0.2251   & 0.2736   &	-0.1237   & 0.0857   &-0.4406\\
         0         &0        &	 0  		&	 -0.0667   &-0.0172   & 0.1418  	&	  0.2158   & 0.2762&    0.2506\\
         0        & 0        & 	0  		 & -0.0787    &0.0360  	& -0.0661  & -0.0605    &0.0366 &   0.0962
         \end{bmatrix}
\end{equation}

\begin{equation}
B^1 =
\begin{bmatrix}
    0.6394 &  -0.3201 \\
    0.8742  & -0.1374 \\
    1.7524   & 0.6158
\end{bmatrix}, \
B^2 =
\begin{bmatrix}
    0.9779  &   0.0399 \\
   -1.1153  & -2.4828 \\
   -0.5500  &  1.1587
\end{bmatrix}, \ 
B^3 =
\begin{bmatrix}
   -1.0263\\
    1.1535\\
   -0.7865
   \end{bmatrix}, \
   B = \begin{bmatrix} B^1 & 0 & 0 \\
   				0 & B^2 & 0 \\
				0 & 0 & B^3
				\end{bmatrix}
				\label{eq:B}
\end{equation}
\end{figure*}

\end{document}